\documentclass{article}

\usepackage{amssymb, amsmath}

\usepackage{tikz-cd}

\newtheorem{theorem}{Theorem}[section]

\begin{document}

\title{A localization theorem for equivariant connective K theory}


\author{Jack Carlisle}

\date{}
\maketitle



\begin{abstract}
We identify Greenlees' $C_n$-equivariant connective K theory spectrum $kU_{C_n}$ as an $RO(C_n)$-graded localization of the actual connective cover of $KU_{C_n}$.
\end{abstract}






\section{Introduction}

For a cyclic group $C_n$, there are two important $C_n$-equivariant forms of connective K theory. One is the connective cover $ku_{C_n}$ of Atiyah's periodic K theory $KU_{C_n}$. As its name suggests, $ku_{C_n}$ is obtained from $KU_{C_n}$ by killing all homotopy Mackey functors of degree $<0$, so that the coefficients of $ku_{C_n}$ are
\[
\underline{\pi}_k(ku_{C_n}) = \begin{cases}
\underline{RU} & k \geq 0 \text{ is even}\\
0 & \text{else}.
\end{cases}
\]
While the definition of this spectrum is quite natural, it fails to satisfy some useful properties. For example, $ku_{C_n}$ does not have thom isomorphisms for all equivariant complex vector bundles, and $ku_{C_n}$ does not satisfy the ``completion theorem", in that the comparison map 
\[
ku^{C_n}_* \to ku^{-*}(B{C_n})
\]
is not completion at the augmentation ideal 
\[
J = \text{ker}\left( \begin{tikzcd} ku^{C_n}_* \ar[r,"\text{res}^{C_n}_e"] & ku_* \end{tikzcd} \right).
\]
For these reasons, $ku_{C_n}$ has been viewed as a pathological ${C_n}$-equivariant form of connective K theory.

In his pioneering work (\cite{Greenlees1}, \cite{Greenlees2}, \cite{Greenlees3}, \cite{Schwede2}), Greenlees' constructs a different ${C_n}$-equivariant form $kU_{C_n}$ of connective K theory via the homotopy pullback square
\[ \begin{tikzcd} 
kU_{C_n} \ar[r] \ar[d] & KU_{C_n} \ar[d] \\
F(E{C_n}_+,ku) \ar[r] & F(E{C_n}_+,KU).
\end{tikzcd} \]
This $C_n$-spectrum is better behaved than $ku_{C_n}$, in that
\begin{itemize}
\item $kU_{C_n}$ is complex oriented, so it has thom isomorphisms for all ${C_n}$-equivariant complex vector bundles,
\item $kU^{C_n}_* \to ku^{-*}(B{C_n})$ is completion at the augmentation ideal $J \subset kU^{C_n}_*$, and
\item  $kU^{C_n}_*$ classifies {\it multiplicative} ${C_n}$-equivariant formal group law, in the sense of \cite{CGK1}. 
\end{itemize}

The purpose of this note is to identify $kU_{C_n}$ as a localization of $ku_{C_n}$. Recall that $kU_{C_n}$ is a complex oriented ${C_n}$-spectrum, meaning that $kU_{C_n}$ has thom isomorphisms for all ${C_n}$-equivariant complex vector bundles. In particular, if we consider an irreducible complex $C_n$-representation $\gamma$ as a vector bundle over a point, we obtain a Thom class $u_\gamma \in kU^{C_n}_{2 - \gamma}$. The class $u_\gamma$ is a unit in $kU^{C_n}_\star$, which means that multiplication by $u_\gamma$ determines an equivalence of ${C_n}$-spectra
\[
\begin{tikzcd} 
kU_{C_n} \ar[r,"\simeq"] & \Sigma^{\gamma-2}kU_{C_n}.
\end{tikzcd} 
\]
It turns out that for any such $\gamma$, the class $u_\gamma$ lifts to  $ku^{C_n}_{\star}$, but is not a unit in $ku^{C_n}_\star$. Our main result is that, by inverting the class $u_\alpha \in ku^{C_n}_\star$ associated to a generator $\alpha$ of the character group 
\[\widehat{C_n} \cong \langle \alpha \mid \alpha^n - 1\rangle,\]
one recovers Greenlees' spectrum $kU_{C_n}$. 
\begin{theorem}\label{maintheorem}
If ${C_n}$ is a cyclic group, then the comparison map
\[
ku_{{C_n}} \to kU_{{C_n}}
\]
induces an equivalence of ${C_n}$-spectra 
\[
ku_{{C_n}}\left[u_\alpha^{-1} \right] \simeq  kU_{{C_n}}.
\]
\end{theorem}

\section{Notation}

We work in the $C_n$-equivariant setting for a fixed cyclic group $C_n$ of order $n \geq 0$. If $V$ is an orthogonal or unitary $C_n$-representation, we write $S(V)$ for the unit sphere in $V$, and $S^V$ for the one point compactification of $V$. If $E_{C_n}$ is a $C_n$ spectrum, we write 
\[
\Sigma^V E_{C_n} = E_{C_n} \wedge S^V 
\]
for the $V$th suspension of $E_{C_n}$. The homotopy groups $\pi_*^H(E_{C_n})$ (as $H$ varies over subgroups of $C_n$) assemble to form a $C_n$ Mackey functor, which we denote $\underline{\pi}_*(-)$. We say $E_{C_n}$ is {\it connective} if $\underline{\pi}_k(E_{C_n}) = 0$ for $k < 0$. More generally, given an integer $m$, we say $E_{C_n}$ is $m$-{\it connective} if $\underline{\pi}_k(E_{C_n}) = 0$ for $k < m$. Every $C_n$ spectrum $E_{C_n}$ has a connective cover $\tau_{\geq 0} E_{C_n} \to E_{C_n}$, which is a connective spectrum such that $\underline{\pi}_k(\tau_{\geq 0}E_{C_n}) \to \underline{\pi}_*(E_{C_n})$ is an isomorphism for $k \geq 0$.

\section{Proof of the theorem}

Recall that $kU_{C_n}$ is defined by the homotopy pullback square
\[ \begin{tikzcd} 
kU_{C_n} \ar[r] \ar[d] & KU_{C_n} \ar[d] \\
F(E{C_n}_+,ku) \ar[r] & F(E{C_n}_+,KU).
\end{tikzcd} \]
The canonical map $ku_{C_n} \to KU_{C_n}$ and the homotopy completion map 
\[ku_{C_n} \to F(E{C_n}_+,ku_{C_n}) \simeq F(E{C_n}_+, ku)\] 
determine a comparison map 
\[ku_{C_n} \to kU_{C_n}.\]
Since $u_\alpha \in ku^{C_n}_\star$ maps to a unit in $kU^{C_n}_\star$, this induces a map $f$ as shown below.

\[ \begin{tikzcd} 
ku_{C_n}  \ar[rr] \ar[dr] & & kU_{C_n}   \\
 & ku_{C_n}\left[ u_\alpha^{-1} \right] \ar[ur,swap,dashed,"f"] 
\end{tikzcd} \]

Our goal is to prove that $f$ is an equivalence. Recall that the localization $ku_{C_n}[u_\alpha^{-1}]$ is defined to be the homotopy colimit of the sequence
\[
\begin{tikzcd} 
ku_{{C_n}} \ar[r,"u_\alpha"] & \Sigma^{\alpha - 2} ku_{{C_n}} \ar[r,"u_\alpha"] & \Sigma^{2(\alpha - 2)} ku_{{C_n}} \ar[r,"u_\alpha"] & \cdots \end{tikzcd} .
\]
Since $ku_{C_n}$ is connective, so is $\Sigma^{m\alpha}ku_{C_n}$ for any $m \geq 0$. This implies that $\Sigma^{m(\alpha-2)} ku_{C_n}$ is $(-2m)$ - connective, and so the composite 
\[\Sigma^{m(\alpha-2)} ku_{C_n}  \to ku_{C_n}[u_{\alpha}^{-1}] \to kU_{C_n}\] 
factors through the $(-2m)$ - connective cover of $kU_{C_n}$ as shown below.

\[ \begin{tikzcd} 
\Sigma^{m(\alpha-2)} ku_{C_n}  \ar[r] \ar[d,swap,dashed,"f_m"] & ku_{C_n}[u_{\alpha}^{-1}]  \ar[d,"f"]  \\
 \tau_{\geq -2m}kU_{C_n} \ar[r] & kU_{C_n}
\end{tikzcd} \] 
Altogether, we have a diagram

\[ \begin{tikzcd} 
ku_{{C_n}} \ar[r,"u_\alpha"] \ar[d,"f_0"] & \Sigma^{\alpha - 2} ku_{{C_n}} \ar[r,"u_\alpha"] \ar[d,"f_1"]& \Sigma^{2(\alpha - 2)} ku_{{C_n}} \ar[r,"u_\alpha"] \ar[d,"f_2"] & \cdots  \\
\tau_{\geq 0}kU_{C_n} \ar[r]  & \tau_{\geq -2} kU_{C_n} \ar[r] & \tau_{\geq -4} kU_{C_n} \ar[r] & \cdots,
\end{tikzcd} \]
which, upon taking homotopy colimits, yields the map $f$. For this reason, if we want to show that $f$ is an equivalence, it suffices to prove that each $f_m$ is an equivalence. We do so by induction on $m$, with the base case $m = 0$ holding by the defining property of the connective cover.

Supposing we have proved that $f_0,\dots,f_m$ are equivalences, we will prove that $f_{m+1}$ is an equivalence. In order to do so, we smash the map 
\[\Sigma^{m\alpha} ku_{C_n} \to \Sigma^{m\alpha}kU_{C_n}\]
 with the cofiber sequence 
 \[S(\alpha)_+ \to S^0 \to S^\alpha.\] 
  Our first claim is that 
 \[\Sigma^{m\alpha} ku_{C_n} \wedge S(\alpha)_+ \to \Sigma^{m\alpha} kU_{C_n} \wedge S(\alpha)_+\]
 is an equivalence. We verify this by observing that $\Sigma^{m\alpha} ku_{C_n} \to \Sigma^{m\alpha} kU_{C_n}$ is an equivalence of non-equivariant spectra, hence it is equivalence after smashing with the free $C_n$ orbit $(C_n)_+$. The existence of the cofiber sequence 
 \[
 (C_n)_+ \to (C_n)_+ \to S(\alpha)_+
 \]
 proves that $\Sigma^{m\alpha} ku_{C_n} \to \Sigma^{m\alpha} kU_{C_n}$ is also an equivalence after smashing with $S(\alpha)_+$, as claimed.

 Using our inductive hypothesis and the 5 lemma, we deduce that the homotopy Mackey functor
 \[\underline{\pi}_*(\Sigma^{m\alpha} ku_{C_n} \wedge S^\alpha) = \underline{\pi}_*(\Sigma^{(m+1)\alpha}ku_{C_n})\]
 is $0$ in negative degrees, and 
 \[\underline{\pi}_*(\Sigma^{(m+1)\alpha} ku_{C_n}) \to \underline{\pi}_*(\Sigma^{(m+1)\alpha}kU_{C_n} )\]
  is an isomorphism in non-negative degrees. In other words, 
\[\Sigma^{(m+1)\alpha} ku_{C_n} \to \Sigma^{(m+1)\alpha} kU_{C_n}\]
 is a connective cover, hence so is the composite 
\[\begin{tikzcd}
\Sigma^{(m+1)\alpha} ku_{C_n} \ar[r] &  \Sigma^{(m+1)\alpha} kU_{C_n} \ar[r,"u_{\alpha}^{-(m+1)}"]  \ar[r,swap,"\simeq"] & \Sigma^{2(m+1)}kU_{C_n}.
\end{tikzcd} \]
This implies that the adjoint map
\[
\Sigma^{(m+1)(\alpha-2)} ku_{C_n} \to kU_{C_n}
\]
is a $(-2(m+1))$-connective cover, and so
\[ \begin{tikzcd} 
\Sigma^{(m+1)(\alpha-2)} ku_{C_n} \ar[r,"f_{m+1}"] &  \tau_{ \geq -2(m+1)} kU_{C_n}
\end{tikzcd}  \]
is an equivalence. This completes our inductive step, and we conclude that
\[
\begin{tikzcd} 
ku_{C_n}[u_{\alpha}^{-1}] \ar[r,"f"] &  kU_{C_n}
\end{tikzcd} 
\]
is an equivalence.



\bibliographystyle{amsplain}

\end{document}